\documentclass[12pt]{amsart}

\input bookman.sty
\boldmath

\usepackage{helvet}

\usepackage{graphics}
\usepackage{amssymb}
\usepackage{amsxtra}
\usepackage{amsmath}
\usepackage{mathrsfs}

\usepackage[arrow, matrix]{xy}
\xyoption{frame}

\theoremstyle{plain}

\newtheorem{theo}{Theorem}[section]
\newtheorem{lemm}[theo]{Lemma}
\newtheorem{prop}[theo]{Proposition}
\newtheorem{coro}[theo]{Corollary}

\theoremstyle{definition}

\newtheorem{defi}[theo]{Definition}
\newtheorem{rema}[theo]{Remark}

\newfont{\rmm}{cmr10 scaled 1000}

\newfont{\itt}{cmsl10 scaled 1000}

\newfont{\rM}{cmr10 scaled 1700}

\newcounter{lemma}[section]

\newcounter{tempcounter}

\newcommand{\lb}{\label}

\newcommand{\rrf}[1]{(\ref{#1})}

\newcommand{\rref}{\rrf}

\newcommand{\Lalf}{\wh\L_{\bar\a}}

\newcommand{\ua}{{\underline\alpha}}

\renewcommand{\a}{\alpha}
\renewcommand{\b}{\beta}
\newcommand{\g}{\gamma}
\renewcommand{\d}{\delta}

\newcommand{\ve}{\varepsilon}

\renewcommand{\t}{\theta}
\renewcommand{\l}{\lambda}

\newcommand{\m}{\mu}

\renewcommand{\o}{\omega}

\renewcommand{\Re}{\text{\rm Re }}
\newcommand{\G}{\Gamma}

\renewcommand{\L}{\Lambda}

\newcommand{\Si}{\Sigma}

\renewcommand{\AA}{{\mathcal A}}

\newcommand{\LL}{{\mathcal L}}

\renewcommand{\SS}{{\mathcal S}}

\newcommand{\cc}{{\mathbb{C}}}

\newcommand{\qq}{{\mathbb{Q}}}
\newcommand{\rr}{{\mathbb{R}}}

\newcommand{\zz}{{\mathbb{Z}}}

\newcommand{\gI}{{\mathfrak{I}}}

\newcommand{\gR}{{\mathfrak{R}}}

\newcommand{\Ker}{\text{\rm Ker }}

\newcommand{\Hom}{\text{\rm Hom}}

\newcommand{\rk}{\text{\rm rk }}

\newcommand{\supp}{\text{\rm supp }}
\newcommand{\Int}{\text{\rm Int }}
\newcommand{\grad}{\text{\rm grad}}

\newcommand{\Id}{\text{\rm Id}}

\newcommand{\bere}{\begin{rema}}
\newcommand{\bede}{\begin{defi}}

\renewcommand{\beth}{\begin{theo}}
\newcommand{\bele}{\begin{lemm}}
\newcommand{\bepr}{\begin{prop}}
\newcommand{\beeq}{\begin{equation}}
\newcommand{\bega}{\begin{gather}}
\newcommand{\begaa}{\begin{gather*}}
\newcommand{\been}{\begin{enumerate}}

\newcommand{\bedee}{\begin{defii}}
\newcommand{\bethh}{\begin{theoo}}
\newcommand{\belee}{\begin{lemmm}}
\newcommand{\beprr}{\begin{propp}}

\newcommand{\beco}{\begin{coro}}

\newcommand{\beal}{\begin{aligned}}

\newcommand{\enre}{\end{rema}}

\newcommand{\enco}{\end{coro}}
\newcommand{\enpr}{\end{prop}}
\newcommand{\enth}{\end{theo}}
\newcommand{\enle}{\end{lemm}}
\newcommand{\enen}{\end{enumerate}}
\newcommand{\enga}{\end{gather}}
\newcommand{\engaa}{\end{gather*}}
\newcommand{\eneq}{\end{equation}}
\newcommand{\enal}{\end{aligned}}

\newcommand{\bq}{\begin{equation}}
\newcommand{\bqq}{\begin{equation*}}

\renewcommand{\leq}{\leqslant}
\renewcommand{\geq}{\geqslant}

\newcommand{\wi}{\widetilde}

\newcommand{\ove}{\overline}

\newcommand{\wh}{\widehat}

\newcommand{\sm}{\setminus}
\newcommand{\ems}{\varnothing}
\newcommand{\sbs}{\subset}

\newcommand{\tens}[1]{\underset{#1}{\otimes}}

\newcommand{\wrt}{with respect to}
\newcommand{\ho}{homomorphism}

\newcommand{\nei}{neighbourhood}

\newcommand{
\sma}{submanifold}

\newcommand{\Prf}{{\it Proof.\quad}}

\newcommand{\smo}{C^{\infty}}

\newcommand{\chart}{\Phi_p:U_p\to B^n(0,r_p)}
\newcommand{\atlas}{\{\Phi_p:U_p\to B^n(0,r_p)\}_{p\in S(f)}}

\newcommand{\qs}{\hfill\square}

\newcommand{\pa}{\vskip0.1in}

\newcommand{\liminv}{\underset {\leftarrow}{\lim}}

\newcommand{\arrh}[3]
{
\xymatrix{
{#1} \ar[r]^<<<<{#2}  &{#3}
}
}

\newcommand{\arrr}[1]
{\arrh {}{#1}{}}

\newcommand{\arrto}
{\xymatrix{{} \ar@{|-{>}}[r]  & {} } }

\newcommand{\arrinto}
{\xymatrix{{} \ar@{^{(}->}[r]  & {} } }

\begin{document}

\title{Circle-valued Morse theory for complex hyperplane arrangements}
\author{Toshitake Kohno and  Andrei Pajitnov}
\address{IPMU, Graduate School of Mathematical Sciences, the
University of Tokyo, 
3-8-1 Komaba, Meguro-ku, Tokyo 153-8914, Japan}
\email{ kohno@ms.u-tokyo.ac.jp}
\address{Laboratoire Math\'ematiques Jean Leray 
UMR 6629,
Universit\'e de Nantes,
Facult\'e des Sciences,
2, rue de la Houssini\`ere,
44072, Nantes, Cedex}                    
\email{ pajitnov@math.univ-nantes.fr}
\begin{abstract}
Let $\AA$ be an essential  complex hyperplane arrangement 
in $\cc^n$, and $H$ denote the union of the hyperplanes.
 We develop the real-valued and circle-valued Morse theory 
on the space $M=\cc^n\sm H$ and prove, in particular, that
$M$ has the homotopy type of a space
obtained from a 
manifold $V$, fibered over a circle,
by attaching to it $|\chi(M)|$ cells of dimension $n$.
We compute the Novikov homology $\wh H_*(M,\xi)$ for 
a large class of homomorphisms $\xi:\pi_1(M)\to\rr$. 
\end{abstract}
\maketitle

\section{Introduction}
\label{s:intro}

Let 
$f$ be a holomorphic Morse
function without zeros
on a complex analytic manifold.
  It gives rise to a real-valued Morse 
function $z\arrto |f(z)|$ and
a circle-valued Morse function $z\arrto f(z)/|f(z)|$.
These two functions can be used to study the topology of 
the underlying manifold. There are, however, numerous technical problems,
and this approach works only in some rare particular cases. This paper
is about one of such cases, namely, the case of the complement to
a complex hyperplane arrangement in $\cc^n$.

Let $\xi_i:\cc^n\to\cc$ be 
non-constant affine functions $(1\leq i \leq m)$;
put $H_i=\Ker\xi_i$. Denote by $\AA$ 
the  hyperplane arangement $\{H_1, \ldots , H_m\}$ and put
$$
H=\bigcup_i H_i, 
\qquad
M(\AA)=\cc^n\sm H.
$$
We will abbreviate $M(\AA)$ to $M$.
The {\it rank} of $\AA$ is the 
maximal codimension of a non-empty 
intersection of some subfamily of $\AA$.
We say that $\AA$ is {\it essential} if $\rk L=n$.
Assume that $\AA$ is essential.
We prove that
$M$ has the homotopy type of a space
obtained from a  finite $n$-dimensional 
CW complex fibered over a circle,
by attaching $|\chi(M)|$ cells of dimension $n$, and 
apply these results to the computation of the Novikov
homology of $M$.

The homology of rank one local systems over $M$ has been
extensively studied. It was shown that this homology
vanishes except in the case of dimension $n$ for generic local systems
(see \cite{AomotoV}, \cite{KohnoHomLoc}).  These homology groups play an 
important role in the theory of hypergeometric integrals and have been
investigated in relation with the cohomology of the twisted de Rham 
cohomology
of logarithmic forms (see \cite{AomotoKita}, \cite{EsnaultViehweg}).

\bede\label{d:positive}
A homomorphism $\xi:\pi_1(M)\to\rr$ is called {\it positive}
if its value on the positive meridian of each hyperplane $H_i$ 
is strictly positive. 
\end{defi}
We show that the structure of the Novikov homology
$\wh H_*(X,\xi)$ for positive homomorphisms $\xi$
is similar to that of the generic local coefficient homology, namely, 
it vanishes in all degrees except 
$n$.

\section{Main results}
\lb{s:results}

Let $\AA$ be an essential arrangement.
Let $\a=(\a_1,\ldots , \a_m)$ be a string of complex numbers.
P. Orlik and H. Terao \cite{OrlikTeraoCr} proved that for 
$\a$ outside a closed algebraic subset of $\cc^m$
the multivalued holomorphic function
$\Phi_\alpha=\xi_1^{\a_1}\cdot \xi_2^{\a_2}\cdot ... \cdot \xi_m^{\a_m}$
has only non-degenerate critical points
(see the works of K. Aomoto 
\cite{AomotoV}, 
and A. Varchenko \cite{Varchenko}
for partial results in this direction).

In this  paper we  work only with $\a\in\rr^m$. It follows
from the Orlik-Terao theorem that there is an open dense
subset $W\sbs \rr^m$ such that for  
$\a\in W$ the function $\Phi_\alpha$ has only 
non-degenerate critical points.
Consider a real-valued  $\smo$ function
$$f_\a(z)=\prod_i \big|\xi_i(z)\big|^{\a_i}, 
\qquad f_\a: \cc^n\sm H \to \rr.$$

\bele\lb{l:morse_absvalue}
Let $\a\in W$.
Then $f_\a$ is a Morse function. The index of every 
critical point of $f_\alpha$
equals $n$.
\enle
\Prf Let $\o_\a=\sum_i \a_i\frac {d\xi_i}{\xi_i}$.
 Then $f_\a(z)=C\cdot \exp(\Re \int_{z_0}^z \o_\a)$, 
therefore, $\log f_\a(z)$ 
locally is the real part of a holomorphic Morse 
function.
In general, if $h$ is a holomorphic Morse function
on an open subset of $\cc^n$, 
then  $\Re h$ 
is a  real-valued Morse function, 
and the index of every critical point of $\Re h$
equals $n$. Our assertion
follows. $\qs$

Let $\ve>0$  and put 
\begin{equation}\lb{f:def_mv}
V
=
\{z\in\cc^n~|~ f_\a(z)=\ve \},\quad
N
=
\{z\in\cc^n~|~ f_\a(z)\geq \ve \}.
\end{equation}
\bede
A vector $\a\in \rr^m$ is called {\it positive}
if $\a_i>0$ for all $i$. The set of all positive vectors
is denoted by $\rr^m_+$.
The {\it rank} of the vector 
$\a\in\rr^m$ is the dimension of the $\qq$-vector space 
generated by the components of $\a$ in $\rr$.
\end{defi}

Recall that we denote $\cc^n\sm H$ by $M$.

\beth\lb{t:hom_type}
Let $\a$ be any positive  vector.
Then for every $\ve>0$ small enough 

1) The inclusion $N\sbs M$
is a homotopy equivalence. 
The space $V=\partial N$ is a $\smo$ manifold of dimension $2n-1$.

2) The space $N$ has the homotopy type
of the space $V$ with $|\chi(M)|$
cells of dimension $n$  attached.

3) If $\a$ has rank $1$, then $V$  is fibered over a circle and 
the fiber has the homotopy type of a finite CW-complex of dimension 
$n-1$.
\enth

To state the next theorem we 
recall  the definition of the Novikov homology.
Let $G$ be a group, and
$\xi:G\to\rr$ a \ho.
Put $G_C=\{g\in G~|~\xi(g)\geq C\}$.
The Novikov completion $\wh\L_\xi$ of the group ring
$\L=\zz G$ \wrt~ the \ho~ $\xi:G\to\rr$ is defined as follows
 (see the
thesis of J.-Cl. Sikorav \cite{SikoravThesis}):
\begin{gather*}
\wh\L_\xi
=\Big\{\l=\sum_{g\in G}n_g g ~\Big|~ \text{ where }~
n_g\in\zz
~ \text{ and }~ \\
\supp\l\cap G_C ~ 
\text{ is finite for every  }~C \Big\}.
\end{gather*}
Let $X$ be a connected topological space and denote 
$\pi_1(X)$ by $G$. Let $\xi:G\to \rr$ be a homomorphism.
The Novikov homology 
$\wh H_*(X,\xi)$
 is by definition 
the homology of the chain complex
$$
\wh\SS_*(\wi X) = \wh\L_\xi \tens{\L}\SS_*(\wi X)
$$
where $\SS_*(\wi X)$ is the singular 
chain complex of the universal
covering of $X$.

\beth\lb{t:nov_vanish}
For any positive homomorphism $\xi:\pi_1(M)\to\rr$ 
the Novikov homology
$\wh H_k(M, \xi)$ vanishes for $k\not=n$ 
and is a free $\wh\L_{\xi}$-module 
of rank $|\chi(M)|$ if $k=n$.
\enth

\section{The gradient field in the \nei~ of $H$}
\label{s:gr_f_neigh}

In this section we develop the main technical 
tools of the paper. The results of this section are valid 
for every arrangement, essential or not.
Let
$$
v_\a(z)
=
\frac {\grad f_\a(z)}{f_\a(z)}.
$$
Denote by $u_j$ the gradient of the function
$z\mapsto |\xi_j(z)|$.
Then
$$
v_\a(z)=
\sum_{j=1}^m 
\a_j
\frac{u_j(z)}{|\xi_j(z)|}.$$
For a linear form 
$\beta:\cc^n\to\cc,\ \  \beta(z)=a_1z_1+\ldots +a_nz_n$ 
the gradient of the function $|\beta(z)|$
is easy to compute, namely
\begin{equation}\lb{f:xi}
\grad|\beta(z)|=\frac {\beta(z)}{|\beta(z)|}\cdot
\big(\bar a_1, \bar a_2,\ldots , \bar a_n \big)
\end{equation}
(it follows, in particular, that 
the norm of this gradient is constant).

\bele\lb{l:grad2}
Assume that the intersection of all 
hyperplanes of $\AA$ is 
non-empty.
Let $\G\sbs\rr_+^m$ be a compact subset.
Then there is $K>0$ such that
\begin{equation}\lb{f:inequal}
||v_\a(z)||\geq K\sum_i \frac 1{|\xi_i(z)|}
\end{equation}
for every $z\in \cc^n\sm H$ and every $\a\in \G$.
\enle
\Prf
We can assume that the intersection of the hyperplanes  
contains $0$, that is, the arrangement is central.
Furthermore, it suffices to prove the Lemma 
for the case when 
\begin{equation}
\lb{f:inter0}
\bigcap_j \Ker\xi_j=\{0\}.
\end{equation}
 Indeed, 
let $L=\cap_j \Ker\xi_j$. Then it follows from 
\rref{f:xi}
that both sides of our inequality \rrf{f:inequal}
are invariant 
\wrt~ translations by vectors in $L$, and 
it is sufficient to prove the formula for
the vector field $v_\a|_{L^\bot}$. 

Furthermore, 
the  both sides of the inequality are 
homogeneous of degree $-1$, and it is sufficient to prove
the inequality for $z\in \Si\sm H$, where $\Si$ stands for 
the sphere of radius $1$ and center $0$.

We will proceed by induction on $m$.
The case $m=1$ being obvious, we will assume that $m>1$.
Choose some $\varkappa>0$, and for $ i\not= j$ let
$$
U_{i,j}
=
\Big\{z\in \Si~\Big|~ |\xi_i(z)|< \varkappa|\xi_j(z)|\Big\}.
$$
These are open sets 
and it follows from the condition
\rrf{f:inter0}
that 
their union $U=\cup_{i,j}U_{i,j}$
covers the set $H\cap \Si$.
We will now prove 
\rrf{f:inequal}
for $z\in U_{i,j}\sm H$. To simplify the notation
let us assume $i=1, j=m$.
Put
\begin{gather*}
A_m(z)= 
v_\a(z), \ \ \ \ 
A_{m-1}(z)=
\sum_{j=1}^{m-1} 
\a_j
\frac{u_j(z)}{|\xi_j(z)|},
\\
B_m(z) = 
\sum_{i=1}^m \frac 1{|\xi_i(z)|},
\
\ \ \ 
B_{m-1}(z) = 
\sum_{i=1}^{m-1} \frac 1{|\xi_i(z)|}.
\end{gather*}
By the induction assumption 
we have $||A_{m-1}(z)||\geq DB_{m-1}(z)$,
where $D$ is some positive constant.
An easy computation shows that
for $z\in U_{i,j}$ we have
$$
||A_m(z)||\geq 
%\frac {D-\varkappa K}{1+\varkappa}
(D-\varkappa\a_m K_m-\varkappa D)B_m(z)
$$
where $K_m=||u_m(z)||$.
Choosing $\varkappa$ sufficiently small,
we conclude that 
$||v_\a(z)||\geq D'B_m(z)$ with  some $D'>0$ for every $\a\in \G$
and $z\in U\sm H$.

The complement $\Si\sm U$ is compact and 
the proof of the Lemma will be over
once we show that $v_\a(z)\not=0$ for $z\in\Si\sm U$.
This is in turn obvious since
$$f_\a(\m z)=
\mu^{\a_1+\ldots +\a_m} f_\a(z)
\text{ \ \ for \ \ } \mu\in\rr_+,
$$
therefore,
$f'_\a(z)\not=0$ for every $z\notin H$,
since all $\a_i$ are positive.
The proof of  Lemma 
\ref{l:grad2} is
now over. $\qs$

For a subset $X\sbs\cc^n$ and $\d>0$ let us denote by
$X(\d)$ the subset of all $z\in\cc^n$ such that $d(z, X)\leq \d$.

\bepr\lb{p:grad1}
Let $\G\sbs\rr_+^m$ be a compact subset.
There is an open \nei~ $U$ of $H$, and numbers $A,B>0$ such that

1) \ \  For some $\d>0$ the set $H(\d)$
is in $U$.

%2) \ \ $||v_\alpha(z)||\geq A$ for every $z\in U\sm H$ and every $\a\in\G$.

2) \ \ 
For every $z\in U\sm H$ and every $\a,\b\in\G$ we have
\begin{equation}\lb{f:aa}
||v_\alpha(z)||\geq A,
\end{equation}
\begin{equation}\lb{f:bb}
||v_\alpha(z)-v_\beta(z)||\leq B\cdot \max_i|\a_i-\b_i|\cdot ||v_\beta(z)||.
\end{equation}
\enpr
\Prf
For a multi-index $I=(i_1,\ldots, i_s)$ let us denote by $H_I$
the intersection of the hyperplanes $H_{i_1}, \ldots ,H_{i_s}$.
Proceeding by induction on $\dim H_I$ 
we will construct for every $I$
with $H_I\not=\ems$
a \nei~ $U_I$ of the subset $H_I$
such that the properties 1) and 2) of the Proposition hold 
if we replace in the formulae $H$ by $H_I$ and $U$ by $U_I$.
Assume that this is done for 
every $H_J$ with $\dim H_J\leq k-1$;
put
\begin{equation}\lb{f:un-1}
U_{k-1}
=
\bigcup_{\dim H_J\leq k-1} U_J.
\end{equation}
Let $I$ be a multi-index
with $\dim H_I=k$. We will construct 
the \nei~ $U_I$. We will assume that $k>0$.
The proof of our assertion for the case 
$\dim H_I=0$ (the initial step of the induction)
 is similar and will be omitted.

We can assume that the multi-index $I$ includes all the values
of $j$ such that $H_I\sbs H_j$. To simplify the notation
let us assume that $I=(1,2,\ldots, r)$.
Write
\begin{equation}\label{f:two_terms}
v_\a(z)
=
\sum_{j=1}^r
\a_j
\frac{u_j(z)}{|\xi_j(z)|}
+
\sum_{j=r+1}^m
\a_j
\frac{u_j(z)}{|\xi_j(z)|}.
\end{equation}
Let $\m>0$, and consider the subset $U'_\mu=H_I(\m)\sm U_{k-1}$.
For $z\in U'_\mu$ the second term of \rrf{f:two_terms}
is bounded uniformly \wrt~ $\a\in\G$
if $\mu>0$ is sufficiently small.
As for the first term, its norm converges to $\infty$ when
$d(z, H_I)\to 0$ as it follows from  Lemma \ref{l:grad2},
applied to the arrangement 
$\{H_1, \ldots , H_r\}$.
An easy computation shows now that for every $\m>0$ sufficiently small
the inequalities 
\rrf{f:aa} and \rrf{f:bb}
hold for $z\in U'_\mu$ and every $\a,\b\in\G$.
Put $U_I=U'_\mu\cup U_{k-1}$.
The properties 1) and 2) for 
$H_I$ and $U_I$ are now easy to deduce.

The \nei~ 
$$U_{n-1}
=
\bigcup_{\dim H_J\leq n-1} U_J.
$$
satisfy the  properties required
in the statement of our Proposition.
 $\qs$

We will also use the normalized gradient 
$$w_\a(z)=\frac {\grad f_\a(z)}{||\grad f_\a(z) ||}.$$
%It is clear that the integral curves  of $w_\a$ are defined 
%on the whole of $\rr$.
Choose 
a \nei~ $U$ of $H$ 
so that
the conclusion of Proposition
\ref{p:grad1}
holds.

\bepr\lb{p:beta_alfa}
Let $\a\in\rr^m_+$. 
There is $D>0$ such that
\begin{equation}\lb{f:w-and-v}
\big\langle v_\a(z), w_\b(z) \big\rangle \geq D
\end{equation}
for every $z\in U\sm H$ and every positive vector $\b$ with $
\max_i|\a_i-\b_i|$ sufficiently small.
\enpr
\Prf
For  $z\in U\sm H$
we have 
$$
\Big|\big\langle
 v_\a(z) -  v_\b(z), 
w_\b(z)
\big\rangle
\Big|
\leq
 B\cdot \max_i|\a_i-\b_i|\cdot ||v_\beta(z)||.
$$
On the other hand
$
\big\langle v_\b(z), 
w_\b(z)
\big\rangle = ||v_\beta(z)||
$
(since these vector fields are collinear and $||w_\beta(z)||=1$), therefore,
$$
\big\langle v_\a(z), w_\b(z)\big\rangle \geq
(1-B\max_i|\a_i-\b_i|)\cdot
||v_\beta(z)||.
$$
If $\a-\b$ is small enough, then  
the right-hand side of the above inequality  is
greater than a positive constant  again by 
Proposition \ref{p:grad1}. $\qs$

\section{The homotopy type of $M$}
\label{s:hom_type}

In this section we prove the first two assertions of
Theorem \ref{t:hom_type}. We fix a positive vector $\a$.
Choose 
a \nei~ $U$ of $H$ 
so that
the conclusion of Proposition
\ref{p:grad1}
holds.
Observe that 
for $\ve>0$ small enough the  set $f_\a^{-1}([0,\ve])$
is in $U$,
therefore, $\ve$ is a regular level of $f_\a$, and $V$ is a 
\sma~ of $M$ of dimension
$2n-1$.
This proves the second part of the assertion 1).

To prove the first part
we use the shift along the flow lines 
of $w_\a$ to construct  
the deformation retraction of 
$M$ onto $N=f_\a^{-1}\big([\ve,\infty[\big)$. 
If $\ve>0$ is sufficiently small, then 
$M\sm N\sbs U$, and 
for every integral curve $\g(t)$ of $w_\a$
starting at a point 
$x\in M\sm N$ we have
$$\frac d{dt} f_\a(\g(t)) = \big\langle \grad f_\a(\g(t)), \g'(t)\big\rangle =
f_\a(\g(t))\cdot ||v_\a(\g(t))||\geq 
A f_a(x)$$
(for every $t$ such that $\g(t)$ is in the set $U$).
Therefore, this trajectory will reach $f_\a^{-1}(\ve)$,
and our deformation retraction is well-defined.

Moving forward to the assertion 2), let us first 
outline the proof. 
Choose $\b\in\rr^m_+$ so that $f_\b:M\to\rr$
is a Morse function, and $\b$ is sufficiently close to $\a$ 
so that the 
property \rref{f:bb} holds.
We are going to apply the Morse theory to
the restriction of $f_\b$ to the manifold $N$ with boundary $V$.
Our setting differs from the classical one (see \cite{MilnorMorse}), 
essentially in 2 points: A) the manifold $N$ is not compact, B) the function $f_\b$ is not constant on the boundary of $N$. 
The main technical tool to deal with these issues is Proposition 
\ref{p:grad1}, which describes the behaviour
of $\grad f_\b$ in a \nei~ of $H$.
Using the classical Morse-theoretic schema
(cf \cite{MilnorMorse}, \S 3), we 
choose a suitable gradient-like vector field $u_\b$ for $f_\b$,
and use the gradient descent along $(-u_\b)$-trajectories to describe the changement of the homotopy type of sublevel sets 
while crossing critical values.
It is technically convenient to replace the couple $(N,V)$ by 
a suitable thickening.
Let $0<\ve_0 < \ve$ and put 
$$
N'=\{z~|~ f_\a(z)\geq \ve_0\},\qquad
L=\{z~|~ \ve_0\leq f_\a(z)\leq \ve\}.
$$
\bele
The inclusion of pairs 
$(N,V)\sbs (N',L)$ 
is a homotopy equivalence.
\enle
\Prf
The shift along the trajectories of the vector field 
$w_\a$ determines a deformation retraction
of $(N',L)$ onto $(N,V)$.
$\qs$

We proceed to the construction of  a
gradient-like vector field for 
 $f_\b$.
We assume that $\beta$ is sufficiently close to $\a$
 so that the property 
\rref{f:w-and-v}
 in Proposition \ref{p:beta_alfa} holds.
For each critical point $c$ of $f_\b$ choose a \nei~ $R_c$
of $c$ such that $\ove{R_c}\cap \ove{U}=\ems$;
let $R$ denote the union of all $R_c$.
Put 
$$K= \{z~|~ f_\a(z)\leq \ve\}. 
%\qquad  K'=\{z~|~ f_\a(z)\leq \ve_0\}.
$$
Using a partition of unity, it is easy to construct a 
$\smo$
vector field $u_\b$ on $\cc^n$ such that
\been\item
$u_\b~|~N$ is a gradient-like vector field for $f_\b$.
\item
$\supp(w_\b-u_\b)\sbs R\cup K$.
\item
$u_\beta~|~(K-L)=0$.
\item
For $z\in L$ we have 
$u_\b(z)=h(f_\b(z))w_\b(z)$ where $h$ is a 
$\smo$ function with the following properties:
\begin{gather*}
   h(t)=0 \quad \text{if}\quad  t\leq \ve_0, \qquad
 h(t)=1 \quad  \text{if}\quad  t\geq \ve, \\
h'(t)>0 \quad \text{for every} \quad  \ve_0<t< \ve.
  \end{gather*}
\enen

Observe that the vector field $u_\b$ is bounded,
therefore, its integral curves  are defined on the whole of $\rr$. 

Now we can investigate the sublevel sets of the function $f_\b$.
For $b\in\rr$ put 
$$
Y_b=\{z~|~ f_\b(z)\leq b\ \  {\rm{ and }} \ \ f_\a(z)\geq \varepsilon_0\}.
$$
It follows from the construction of the vector field $u_\b$, that 
$Y_b$ is $(-u_\b)$-invariant, that is, 
every
$(-u_\b)$-trajectory $\g(z,t)$  starting at 
a point $z\in Y_b$ remains in $Y_b$ for all $t\geq 0$.
Furthermore, it follows from Proposition 
\ref{p:beta_alfa} that $L$ is also 
$(-u_\b)$-invariant. Therefore, $Y_b\cup L$ is $(-u_\b)$-invariant.

\bepr\label{p:nocrit}
Let  $0\leq a<b$ and assume that $f_\b$ does 
not have critical points in $Y_b\sm Y_a$.
Then 
\been\item
There is $C>0$ such that for every $z\in Y_b\sm (Y_a\cup L)$
we have 
\begin{equation}\lb{f:greater}
f'_\b(z)(u_\b(z)) \geq C.
\end{equation}
\item
There is $T>0$ such that 
$\g(z,T)\in \Int L \cup \Int Y_a$
for every $z\in Y_b$.
\item
The inclusion
$$(Y_a\cup L, L)\sbs (Y_b\cup L, L)$$
is a homotopy equivalence.
\enen
\enpr
\Prf
1) \quad
Consider two cases: A) $z\notin U$, and B) $z\in U$.
For the case A) observe that the 
set $Y_b\sm (\Int Y_a \cup U)$ is compact
(since the arrangement $\AA$ is essential),
the function $f_\b$ has no critical points in it
and $u_\b$ is a gradient-like vector field for 
$f_\b$. Thus the property \rref{f:greater}
is true for
$z\in Y_b\sm (\Int Y_a \cup U)$.

For the case B) observe that when 
$z\in (Y_b\cap U)\sm (Y_a\cup L)$
we have $u_\b(z)=w_\b(z)$ and 
$$
f'_\b(z)(u_\b(z))
=
\langle
v_\b(z)f_\b(z), w_\b(z)
\rangle
\geq
a||v_\b(z)|| \geq aA$$
(see Proposition \ref{p:grad1}).
\pa
2) Let $T>\frac {b-a}C$, and 
assume that 
$\g(z, T)\notin \Int L \cup \Int Y_a$.
The set $Y_a$ is obviously invariant
with respect to the shift along $(-u_\b)$-trajectories.
The same is true for $L$, as it follows from 
\rref{f:w-and-v}.
Therefore,
$\g(z, \tau)\notin \Int L \cup \Int Y_a$
for all $\tau\in [0,T]$.
By the first assertion of the Lemma we have 
$$
f_b(\g(z, T))<
f_b(\g(z, 0))
-C\cdot \frac {b-a}C
\leq a
$$
which contradicts the assumption 
$\g(z, T)\notin \Int Y_a$.
\pa
3) Define a homotopy 
$$H_t:Y_b\cup L \to Y_b\cup L,\qquad t\in [0, T],\qquad
H_t(z)=\g(z,t).$$
The sets 
$Y_b, L, Y_a$ are $(-u_\b)$-invariant 
and it follows from the second assertion of the Lemma that $H_T(z)\in Y_a\cup L$
for every $z\in Y_b\cup L$.
Applying the next Lemma, we 
accomplish the proof of the Proposition.

$\qs$

\bele\lb{l:hom_eq}
Let $(X,A)\sbs (Y,B)$
be pairs of topological spaces, and $H_t:Y\to Y$
be a homotopy, where $t\in [0, T]$ such that
\been\item
$X, A, B$ are invariant under $H_t$ for each $t$.
\item $H_0=\Id, \ H_T(Y)\sbs X,\ H_T(B)\sbs A$.
\enen
Then the inclusion $(X,A)\arrinto (Y,B)$
is a homotopy equivalence.
\enle
\Prf The proof repeats the arguments of 
\cite{PajitnovCVMT}, 
Lemma 1.8, page 171. $\qs$
\beco
Assume that $f_\b$ does not have 
critical points in $Y_a$. Then the inclusion
$L\sbs Y_a\cup L$ is a homotopy equivalence. $\qs$
\enco

\pa

Now let $c$ be a critical value of $f_\b$;
we will describe the homotopy type of the pair 
$(Y_{c+\d}\cup L, L)$
in terms of $(Y_{c- \d}\cup L, L)$.
Denote 
by $p_1, ..., p_k$ the critical points on the level $c$
and let $H_1, ... , H_k$ be the corresponding 
$(-u_\b)$-invariant 
handles
around the critical points; denote 
the union $\cup_i H_i$ by $H$.
If $\d>0$ is sufficiently small, then 
we can choose the handles in such a way that 
they  do not intersect $L$, and, therefore, the set
$$
L\cup Y_{c-\d}\cup H
$$
is homeomorphic to the result 
of attaching of $k$ handles of index $n$
to $L\cup Y_{c-\d}$.
%Let $H=\cup_i H_i$.

\bele\label{l:crit}
\been\item
There is $C>0$ such that for every
$$z\in Y_{c+\d}
\sm (Y_{c-\d}\cup H \cup L )$$
we have
$f'_\b(z)(u_\b(z))\geq C$.
\item 
There is $T>0$ such that
$$
\g(z,T)\in  \Int \big(Y_{c-\d} \cup  H   \cup L\big)$$
for every $z\in Y_{c+\d}$.
\item
The inclusion 
$$
(Y_{c-\d} \cup H \cup L, L)
\sbs
(Y_{c+\d}\cup L, L)
$$
is a homotopy equivalence.
\enen
\enle
\Prf
The proof repeats the corresponding steps
of the proof of Proposition \ref{p:nocrit}
and will be omitted. $\qs$

The point 2) of Theorem \ref{t:hom_type} follows by the usual 
inductive procedure of Morse theory. Namely,
we recover
 the homotopy type of a manifold by successive 
attaching of handles corresponding to
the critical points of a Morse function on it.

\section{The Novikov homology  of $V$ and $M$}
\label{s:homology}

In this section we prove 
the assertion 3) of Theorem \ref{t:hom_type}
and Theorem \ref{t:nov_vanish}.

Returning to the space $M=\cc^n\sm H$, observe that
$H_1(M,\zz)$ is a free abelian group of rank
 $m$ generated by the meridians of the hyperplanes $H_i$.
The  elements of the dual basis in the 
group $\Hom (\pi_1(M),\zz)$
will be denoted by $\t_i$ (where $1\leq i \leq m$), 
%so that
%$\t_i$ is represented by the differential 
%form $\frac 1{2\pi\sqrt{-1}}\frac {d\xi_i}{\xi_i}$. 
For $\a\in\rr^m$ 
denote by $\ove \alpha:\pi_1(M)\to\rr$ the \ho~ 
$\sum_i\a_i\t_i$.
It is clear that $\alpha\in\rr^m$ is positive if and only if
$\ove \alpha$ is a  positive homomorphism.
The composition 
$$\pi_1(N)\arrr {\approx}\pi_1(M) \arrr {\ove \alpha} \rr,$$
where the first homomorphism is induced by the inclusion 
$N\sbs M$,
will be denoted by the same symbol $\ove \alpha$
by an abuse of notation. We will restrict ourselves to the case 
$n\geq 2$, where $V$ is connected. The proof in the case $n=1$ is easy 
and will be omitted. 
The composition 
$$\pi_1(V)\to\pi_1(M) \arrr {\ove \alpha} \rr,$$
where the first homomorphism is induced by the inclusion 
$V\sbs M$,
will be denoted by $\ua$.
Recall the holomorphic 1-form 
$$\o_\a=\frac{d\Phi_\alpha}{\Phi_\alpha} = 
\sum_j \a_j\frac {d\xi_j}{\xi_j}$$
and denote its real and imaginary parts by
 $\gR$ and $\gI$ respectively.
Then $\gR=\frac {df_\a}{f_\a}$.
Denote by $\wh\alpha$ the element
of the group $H^1(M,\rr)$ corresponding to $\ove\alpha$
under the canonical isomorphism
$$\Hom (\pi_1(M),\rr)\approx H^1(M,\rr).$$
Then 
the cohomology 
class of $\gI$ equals $2\pi \wh\a$.
Let $\iota_\a$ be the vector field dual to $\gI$.
Since $\o_\a$ is a holomorphic form,
we have
\begin{equation}\lb{f:real_im}
\frac {\grad f_\a(z)}{f_\a(z)}
= -\sqrt{-1} \cdot \iota_\a(z).
\end{equation}
If $\ve>0$ is small enough 
so that $V$ is contained in the \nei~ $U$ 
from Proposition
\ref{p:grad1} we have 
$||\iota_\a(z)||\geq A$. 
Observe that  $\iota_\a(z)$ is 
orthogonal to $\grad f_\a(z)$
and, therefore, tangent to $V=f_\a^{-1}(\ve)$.
 We deduce
that the restriction to $V$ of the 1-form $\gI$
does not vanish, and, moreover, the norm
 of its dual vector field
is bounded from below by a strictly positive constant.

\bepr\lb{p:novhom_V}
The Novikov homology 
$\wh H_k(V,\ua)$ 
vanishes for all $k$.
\enpr
\Prf
Let $p:\wi V\to V$ be the universal covering of $V$. 
The closed $1$-form $p^*(\gI)$ is cohomologous to zero;
let $p^*(\gI)=dF$, where 
$$F=
{\rm{arg}}\Big(
\xi_1^{\a_1}\cdot \xi_2^{\a_2}\cdot ... \cdot \xi_m^{\a_m}\Big)
$$
is a real-valued function on $\wi V$
without critical points. Denote by $B_n$
the subset $F^{-1}\big(]-\infty, -n]\big)$.
The singular chain complexes $C_*^{(n)}=\SS_*(\wi V, B_n)$
form an inverse system, and 
the Novikov homology
$\wh H_*(V,\ua)$
is isomorphic to the homology of its inverse limit 
%$\liminv\ C_*^{(n)}$
(see \cite{SikoravThesis}).
For every $k$ we have an exact sequence
$$
{\lim}^1 H_{k+1}\big(C_*^{(n)}\big) \to  H_k\big(\liminv\ C_*^{(n)}\big)
\to 
\liminv\ H_k\big(C_*^{(n)}\big).
$$
The lift of the vector field $\iota_\a$
to $\wi V$ will be denoted by the same letter
$\iota_\a$.
The  standard argument using the shift diffeomorphism 
along the trajectories  of $-\iota_\a$
shows that $H_k(C_*^{(n)})=0$ for every $k$; 
our Proposition follows. $\qs$

Consider now the case when $\a$ is of rank one, that is, 
all $\a_i$ are rational multiples of one real number.
In this case the differential form 
$\gI$ 
is the differential of a map 
$g:V\to \rr/a\rr\approx S^1$ for some $a>0$.

\bepr\lb{p:fibr}
The map $g$ is a fibration of \ $V$ over $S^1$.
\enpr
\Prf
The map $g$ does not have critical points.
Consider the vector field
$$
y_\a(z)= 
\frac {\iota_\a(z)}{||\iota_\a(z)||^2}.$$
For $x\in V$ denote by $\g(x,t;y_\a)$
the $y_\a$-trajectory starting at $x$.
Since the norm of  $\iota_\a(z)$
is bounded away from zero in $V$,
the trajectory is defined on the whole of $\rr$.
We have also
$$
\frac d{dt} \Big(g\big(\g(x,t;y_\a)\big)\Big)=1.
$$
Pick any 
$b\in\rr$ and  let $\bar b$ be the image of $b$ under the projection
$\rr\to \rr/a\rr \approx S^1$; let $c\in S^1$ be the image of
$b+a/2$. Let $V_0=g^{-1}(\bar b)$. It is easy to check
that the map 
$$
(x,t)\mapsto 
\g(x,t;y_\a)
$$
is a diffeomorphism 
$$V_0\times \big]-a/2,a/2\big[  
\approx\ 
g^{-1}\big(S^1\sm \{ c\}\big)$$
compatible with projections.
Therefore, $g$ is a locally trivial fibration.
$\qs$

It is clear that  any fiber of $g$ is 
locally the set of zeros of a holomorphic 
function, therefore, it is a closed complex analytic submanifold of
$\cc^n$ and  has a homotopy type of a 
CW-complex of dimension $\leq n-1$
(see \cite{AndreottiFrankelLefHyp}). 
Moreover, it is not difficult to prove 
that for any $a\in\cc^n$ the distance function $d(z)=||z-a||^2$
restricted to $X$ has no critical points outside
a ball $D(a,R)$ of sufficiently large raduis $R$. 
Therefore,  the Morse theory applied to the function $d|X$, implies
that $X$ has the homotopy type of a finite CW-complex 
of dimension $\leq n-1$. A similar argument  shows that 
the manifold $V$ itself is homotopy equivalent to a finite CW complex.
The proof of Theorem \ref{t:hom_type} is now complete. $\qs$

{\it Proof of Theorem \ref{t:nov_vanish}}. \
Let
%by $G$ the fundamental group of $N$, 
$\L$  be the  group ring of the fundamental group of $N$.
Let $\xi:\pi_1(M)\to\rr$ be a posiive homomorphism, then $\xi=\ove\alpha$
where $\alpha\in\rr^m$ is a positive vector.
Let
$\Lalf$ denote the Novikov 
completion of $\L$ \wrt~ $\ove\a$.
Denote by $q:\wi N\to N$ the universal covering of $N$.
Put $\bar V=q^{-1}(V)$.

We have the short exact sequence of free $\Lalf$-complexes
\begin{equation}\lb{f:sh_ex}
0
\to 
\Lalf\tens{\L} 
\SS_*(\bar V)
\to 
\Lalf\tens{\L}
\SS_*(\wi N)
\to 
\Lalf\tens{\L}
\SS_*(\wi N, \bar V)
\to 
0
\end{equation}
Observe that the inclusion $V\sbs N$ induces a surjective homomorphism 
of fundamental groups. Therefore, the space $\bar V$ is connected
and the covering $\bar V\to V$ is a quotient of the 
universal covering $\wi V\to V$,
so that 
$$ H_*\big(\Lalf\tens{\L}\SS_*(\bar V)\big) = 
H_*\Big(\Lalf\tens{\wh\LL_\ua}  
\big( {\wh\LL_\ua}\tens{\LL}\SS_*(\wi V)\big)\Big),$$
where 
$\LL$ is the group ring of the fundamental group of $V$,
and $\wh\LL_\ua$ is the Novikov completion of $\LL$
\wrt~ $\ua$.
By Proposition \ref{p:novhom_V} 
the Novikov homology of $V$ 
vanishes. Consider  the  long  exact 
sequence of homology modules, 
derived from the short exact sequence
\rrf{f:sh_ex}. 
Since 
$H_*\big(\Lalf\tens{\L} 
\SS_*(\bar V)\big)=0$, we deduce 
the isomorphisms
$$
\wh H_*(N, \bar\a)
\approx
H_*\Big(\Lalf\tens{\L}
\SS_*(\wi N, \bar V)\Big).
$$
Observe now that the homology of the couple $(\wi N, \bar V)$
is a free module over $\L$ of rank $|\chi(M)|$,
concentrated  in degree $n$. 
Theorem \ref{t:nov_vanish} follows. $\qs$

\section{Non-essential arrangements}
\label{s:further}

Let us consider the case of  non-essential arrangements $\AA$.
%As before we denote by $m$ the number of hyperplanes in the 
%family $\AA$.
Assume that  $\rk \AA=l<n$. 
The function $f_\a$ is not a Morse function in this case.
However, the analog of Theorem \ref{t:hom_type} is easily obtained by reduction to 
the case of essential arrangements.

Denote by $\pi:\cc^l\oplus \cc^k\to \cc^l$ the projection onto 
the first direct summand.
Let $\AA$ be a hyperplane arrangement in $\cc^l$, defined by affine
functions $\xi_i:\cc^l\to\cc$.
The functions $\xi_i\circ \pi$ determine a hyperplane 
 arrangement in $\cc^{l+k}$
which will be called {\it $k$-suspension } of $\AA$. 
It is not difficult to prove the next proposition.

\bepr\lb{p:reduc}
Let $\AA$ be a hyperplane arrangement in $\cc^n$ of rank $l$.
The $\AA$ is linearly isomorphic to the $(n-l)$-suspension 
of an essential hyperplane arrangement $\AA_0$ in $\cc^l$.
\enpr

The complement 
$M(\AA)$ 
is  diffeomorphic 
to
$M(\AA_0)\times\cc^{n-l}$. Denote $M(\AA)$ by $M$ for brevity.
We obtain the following 
generalizations of  theorems
\ref{t:hom_type}
and
\ref{t:nov_vanish}.

\beth\lb{t:hom_type_noness}
There is an $l-1$-dimensional manifold $Y$,\
fibered over a circle, such that 
$M$ 
is homotopy equivalent to the
result of attaching  to $Y$ of 
$|\chi(M)|$ cells of dimension
$l$.
\enth

\beth\lb{c:novhom_noness}
For every positive homomorphism $\xi:\pi_1(M)\to \rr$
the Novikov homology
$\wh H_k(M, \xi)$ vanishes for every $k\not=l$ and
is a free module of rank 
$|\chi(M)|$ for $k=l$.
\enth

\begin{rema}
Theorem \ref{t:hom_type}
leads to a quick proof of the  well-known 
fact that the space $M$ is homotopy equivalent to
a finite CW-complex. 
\end{rema}

\begin{rema}
A. Suciu communicated to us that our results
are related to the recent work 
\cite{PapadimaSuciuBNSjump} of 
S. Papadima and A. Suciu.
Let $A$ denote the Orlik-Solomon algebra for 
the arrangement $\AA$, which is isomorphic to 
the cohomology ring of $M(\AA)$. We denote by 
$A^j$ the degree $j$ part of $A$. Now the Aomoto
complex is the cochain complex
$A^* = \oplus_{j \geq 0} A^j$ whose coboundary operator
is defined to be the multiplication by $a \in A^1$. 
The corresponding 
cohomology is denoted by $H^*(A, a)$.
Let us call an element $a$ {\it resonant}
in degree $j$ 
if the cohomology 
$H^j(A, a)$ does not vanish.
Let $b\in\rr^m$;
recall from Section 
%%%%%%
\ref{s:homology}
the corresponding homomorphism
$\ove b:\pi_1(M)\to\rr$ 
and the cohomology class
$\wh b=\sum_i b_i \theta_i\in H^1(M,\rr)$.

Proposition 16.7 of the paper \cite{PapadimaSuciuBNSjump}
says, that  if 
the Novikov 
homology 
$\wh H_j(M(\AA), \ove{b})$
equals $0$ 
for $j \leq q$, then
 $\wh b$ is non-resonant
in degrees $j \leq q$.
Thus our Theorem \ref{t:nov_vanish}
implies that every positive vector $b$ 
the corresponding cohomology class
$\wh b $ 
is non-resonant
in degrees strictly less than the rank of 
the arrangement.
This also follows from a result due to S.~Yuzvinsky
(see \cite{Yuzvinsky1}, \cite{Yuzvinsky2}).
\end{rema}

\pa
\centerline{\bf Acknowledgments}
\pa

The work was supported by World Premier International Research Center
Initiative (WPI Program), MEXT, Japan.

The work was accomplished during  the second author's
stay at the Graduate School of Mathematical Sciences, 
the University of Tokyo in the autumn  of 2010 and 2011.
The second author is grateful 
to the Graduate School of Mathematical Sciences
for warm hospitality.

The authors thank A. Suciu for helpful discussions.


\begin{thebibliography}{99}

\bibitem{AndreottiFrankelLefHyp}
{A.~Andreotti and T.~Frankel,
\emph{
The Lefschetz theorem on hyperplane sections},
Ann. of Math. (3)
{\bf 69} (1959), 713-717}.

\bibitem{AomotoV}
{K.~Aomoto,
\emph{
On vanishing of cohomology attached
to certain many valued meromorphic functions},
J. Math. Soc. Japan (2) {\bf 27} (1975),
248-255}.



\bibitem{AomotoKita}
{K.~Aomoto, M.~Kita, P~.Orlik and H.~Terao, 
\emph{
Twisted de Rham cohomology groups of logarithmic forms},
Adv. Math.  {\bf 128} (1997),
119-152}.



\bibitem{EsnaultViehweg}
{E.~Esnault, V.~Schechtman and E.~Viehweg, 
\emph{
Cohomology of local systems of the complement of hyperplanes},
Invent. Math.  {\bf 109} (1992),
557-561}.


\bibitem{KohnoHomLoc}
{T.~Kohno,
\emph{Homology of a local system
on the complement of hyperplanes},
Proc. Japan Acad. {\bf 62} Ser. A (1986), 144--147}.



\bibitem{MilnorMorse}
{J.~ Milnor,
{Morse theory},
Annals of Mathematics Studies Number 51,
Princeton University Press, 
1963}.



\bibitem{OrlikTeraoBookArr}
{P.~Orlik and H.~Terao,
{Arrangements of Hyperplanes},
Grundlehren der mathematischen
Wissenschaften {\bf 300},
Springer-Verlag
1992}.

\bibitem{OrlikTeraoCr}
{P. Orlik and H. Terao,
\emph{The number of critical points of a product of
powers of linear functions},
Invent. Math.
{\bf 120}
(1995),
1 -- 14}.

\bibitem{PajitnovCVMT}
{A.~V.~Pajitnov,
{Circle-valued Morse Theory},
de Gruyter Studies in Mathematics {\bf 32},
Walter de Gruyter
2006}.

\bibitem{PapadimaSuciuBNSjump}
{S. Papadima, A. Suciu,
\emph{ Bieri-Neumann-Strebel-Renz invariants, and
         homology jumping loci,}
Proc. London Math. Soc. (3) 100 (2010) 795-834.}

\bibitem{SikoravThesis}
{J.-Cl.~Sikorav,
\emph{
Points fixes de diff\'eomorphismes
symplectiques, intersections de sous-vari\'et\'es
lagrangiennes, et singularit\'es de un-formes ferm\'ees},
Th\`ese de Doctorat d'Etat Es
Sciences Math\'ematiques,
Universit\'e Paris-Sud, Centre d'Orsay, 1987}.


\bibitem{Varchenko}
{A.~Varchenko, 
\emph{Critical points of the product of powers of
linear functions and families of bases of singular vectors},
Compos. Math. {\bf 97} (1995), 385-401}.

\bibitem{Yuzvinsky1}
{S.~Yuzvinsky, 
\emph{Cohomology of the Brieskorn-Orlik-Solomon algebras}, 
Comm. Algebra (14) {\bf 23}  (1995), 5339--5354}. 

\bibitem{Yuzvinsky2}
{S.~Yuzvinsky, 
\emph{Orlik-Solomon algebras in algebra and topology}, 
Uspehi Mat. Nauk {\bf 56}  (2001), 87-166}.

\end{thebibliography}
\end{document}